\documentclass[12pt]{report}
\usepackage{amssymb,amsmath,makeidx,verbatim}
\newcommand{\R}{\mathbb{R}}

\newcommand{\C}{\mathbb{C}}
\newcommand{\Z}{\mathbb{Z}}

\newcommand{\be}{\begin{enumerate}}
\newcommand{\ee}{\end{enumerate}}
\newcommand{\bq}{\begin{eqnarray*}}
\newcommand{\eq}{\end{eqnarray*}}

\begin{document}
\newcommand{\disp}{\displaystyle}
\thispagestyle{empty}
\begin{center}
\textsc{The full Bochner theorem on real reductive groups\\}
\ \\
\textsc{Olufemi O. Oyadare}\\
\ \\
Department of Mathematics,\\
Obafemi Awolowo University,\\
Ile-Ife, $220005,$ NIGERIA.\\
\text{E-mail: \textit{femi\_oya@yahoo.com}}\\

\end{center}
\begin{quote}
{\bf Abstract.} {\it The major results of Barker $[3.],$ leading to the spherical Bochner theorem and its (spherical) extension, were made possible through the spherical transform theory of Trombi-Varadarajan $[14.]$ and were greatly controlled by the non-availability of the full (non-spherical) Harish-Chandra Fourier transform theory on a general connected semisimple Lie group, $G.$ Sequel to the recently announced results of Oyadare $[13.],$ where the full image of the Schwartz-type algebras, $\mathcal{C}^{p}(G),$ under the full Fourier transform is computed to be $\mathcal{C}^{p}(\widehat{G}):=\{(\widehat{\xi_{1}})^{-1}\cdot h\cdot (\widehat{\xi_{1}})^{-1}:h\in\bar{\mathcal{Z}}({\mathfrak{F}}^{\epsilon})\}$ with $\bar{\mathcal{Z}}({\mathfrak{F}}^{\epsilon})$ given as the Trombi-Varadarajan image of $\mathcal{C}^{p}(G//K)$ ($0< p\leq2$), the present paper now gives the full Bochner theorem for $G$ by lifting the earlier results of Barker $[3.]$ to full non-spherical status. An extension of the full Bochner theorem to all of $\mathcal{C}^{p}(G),$ $1\leq p\leq2,$ is then established. It is also conjectured that every positive-definite distribution $T$ on $G$ which corresponds to a Bochner measure $\mu$ on ${\mathfrak{F}}^{\epsilon}$ extends uniquely to an element of $\mathcal{C}^{p}(G)'$ if and only if $T$ can be expressed as a finite sum of derivatives of a class of functions exclusively parameterized by members of ${\mathfrak{F}}^{\epsilon}$ and $supp\; (\mu)\subset{\mathfrak{F}}^{\epsilon}$ (with $\epsilon=(\frac{2}{p})-1$ for all $1\leq p\leq2).$ This gives the non-spherical abstract version of the extension theorem for any positive-definite distribution on $G.$ Our results confirm the one-to-one correspondence between tempered invariant positive-definite distributions and the Bochner measures of the case $SU(1,1)/\{\pm1\}$ (as earlier computed by Barker $[5.]$) for all $G.$}
\end{quote}
\ \\

{\bf $\bf{\S 1.}\;$ Introduction.}

The classical \textit{Bochner theorem} is one of the cornerstones of harmonic analysis on $\R^{n}$ and it states that every \textit{positive-definite distribution} on $\R^{n}$ is the \textit{Fourier transform} of some unique \textit{tempered measure} on $\R^{n}.$ The proof of this theorem depends on being able to express such a distribution as a finite sum of derivatives of bounded functions on $\R^{n}$ and showing that the distribution is tempered. Underlying this proof is the use of the Fourier transform theory of \textit{Schwartz space} on $\R^{n}$ which makes it possible to, in the first instance, define a distribution on $\R^{n}.$

In generalizing the classical Bochner theorem to other \textit{topological groups} Barker $[4.]$ showed that every positive-definite distribution on a \textit{unimodular Lie group} may be expressed as a finite sum of derivatives of bounded functions. However, the Fourier transform theory for Schwartz space is not available for unimodular Lie groups neither had the theory been completely successful for specific classes of unimodular Lie groups, except for isolated examples like $SL(2,\R).$ The most successful Fourier transform theory for a class of unimodular Lie groups had been the \textit{Trombi-Varadarajan theory} $[14.]$ constructed only on the \textit{spherical} part of a connected semisimple Lie group $G$ and this informed the consideration of the \textit{spherical Bochner theorem} (which still lacks the classical temperedness) and its \textit{spherical extension} in $[3.]$ (which may be consulted for other partial versions and references).

In this paper, we shall however employ the recently announced results of Oyadare $[13.]$ on the full transform theory for $G$ to give the \textit{full Bochner theorem.} The paper is organized as follows. Section $2$ contains preliminary matters on the \textit{structure theory} of $G,$ explicit parametrization of its \textit{spherical functions} and the \textit{Helgason-Johnson theorem.} The results leading to the spherical Bochner theorem and spherical extension are discussed in section $3,$ where the \textit{Fundamental Theorem of Harmonic Analysis on $G$} (which generalizes the Trombi-Varadarajan theorem) is stated as Theorem $3.2$ $(cf.\;[13.]).$ The last section contain the main contributions of this paper, where the existence of the Fundamental Theorem paves way to the uplifting of Barker's spherical results to all of $G$ (Theorems $4.2$ and $4.3$). An \textit{abstract extension} (Conjecture $4.5$) of our results suggests that the starting point of expressing a positive-definite distribution as a finite sum of derivatives of bounded functions may need to be critically looked at in order to recover their classical temperedness.
\ \\

{\bf $\bf{\S 2.}\;$ Structure of the group $G$ and its spherical functions.}

Let $G$ be a connected semisimple Lie group with finite center, we denote its Lie algebra by $\mathfrak{g}$
whose \textit{Cartan decomposition} is given as $\mathfrak{g} = \mathfrak{t}\oplus\mathfrak{p}.$ Denote by $\theta$ the \textit{Cartan involution} on $\mathfrak{g}$ whose collection of fixed points is $\mathfrak{t}.$
We also denote by $K$ the analytic subgroup of $G$ with Lie
algebra $\mathfrak{t}.$  $K$ is then a maximal compact subgroup of $G.$
Choose a maximal abelian subspace  $\mathfrak{a}$ of $\mathfrak{p}$ with algebraic
dual $\mathfrak{a}^*$ and set $A =\exp \mathfrak{a}.$  For every $\lambda \in \mathfrak{a}^*$ put
$$\mathfrak{g}_{\lambda} = \{X \in \mathfrak{g}: [H, X] =
\lambda(H)X, \forall  H \in \mathfrak{a}\},$$ and call $\lambda$ a restricted
root of $(\mathfrak{g},\mathfrak{a})$ whenever $\mathfrak{g}_{\lambda}\neq\{0\}$.
Denote by $\mathfrak{a}'$ the open subset of $\mathfrak{a}$
where all restricted roots are $\neq 0,$ and call its connected
components the \textit{Weyl chambers.}  Let $\mathfrak{a}^+$ be one of the Weyl
chambers, define the restricted root $\lambda$ positive whenever it
is positive on $\mathfrak{a}^+$ and denote by $\triangle^+$ the set of all
restricted positive roots. Members of $\triangle^+$ which form a basis for $\triangle$ and can not be written as a linear combination of other members of $\triangle^+$ are called \textit{simple.} We then have the \textit{Iwasawa
decomposition} $G = KAN$, where $N$ is the analytic subgroup of $G$
corresponding to $\mathfrak{n} = \sum_{\lambda \in \triangle^+} \mathfrak{g}_{\lambda}$,
and the \textit{polar decomposition} $G = K\cdot
cl(A^+)\cdot K,$ with $A^+ = \exp \mathfrak{a}^+,$ and $cl(A^{+})$ denoting the closure of $A^{+}.$

If we set $M = \{k \in K: Ad(k)H = H$, $H\in \mathfrak{a}\}$ and $M' = \{k
\in K : Ad(k)\mathfrak{a} \subset \mathfrak{a}\}$ and call them the
\textit{centralizer} and \textit{normalizer} of $\mathfrak{a}$ in $K,$ respectively, then (see $[9.]$, p. $284$);
(i) $M$ and $M'$ are compact and have the same Lie algebra and
(ii) the factor  $\mathfrak{w} = M'/M$ is a finite group called the \textit{Weyl
group}. $\mathfrak{w}$ acts on $\mathfrak{a}^*_{\C}$ as a group of linear
transformations by the requirement $$(s\lambda)(H) =
\lambda(s^{-1}H),$$ $H \in \mathfrak{a}$, $s \in \mathfrak{w}$, $\lambda \in
\mathfrak{a}^*_\mathbb{\C}$, the complexification of $\mathfrak{a}^*$.  We then have the
\textit{Bruhat decomposition} $$G = \bigsqcup_{s\in \mathfrak{w}} B m_sB$$ where
$B = MAN$ is a closed subgroup of $G$ and $m_s \in M'$ is the
representative of $s$ (i.e., $s = m_sM$). The Weyl group invariant members of a space shall be denoted by the superscript $^{\mathfrak{w}}.$

Some of the most important functions on $G$ are the \textit{spherical
functions} which we now discuss as follows.  A non-zero continuous
function $\varphi$ on $G$ shall be called a \textit{(zonal) spherical
function} whenever $\varphi(e)=1,$ $\varphi \in C(G//K):=\{g\in
C(G)$: $g(k_1 x k_2) = g(x)$, $k_1,k_2 \in K$, $x \in G\}$ and $f*\varphi
= (f*\varphi)(e)\cdot \varphi$ for every $f \in C_c(G//K),$ where $(f \ast g)(x):=\int_{G}f(y)g(y^{-1}x)dy$ $[7.].$  This
leads to the existence of a homomorphism $\lambda :
C_c(G//K)\rightarrow \C$ given as $\lambda(f) = (f*\varphi)(e)$.
This definition is equivalent to the satisfaction of the functional relation $$\int_K\varphi(xky)dk = \varphi(x)\varphi(y),\;\;\;x,y\in G,\;[10.].$$

It has been shown by Harish-Chandra [$8.$] that spherical functions on $G$
can be parametrized by members of $\mathfrak{a}^*_{\C}.$  Indeed every
spherical function on $G$ is of the form $$\varphi_{\lambda}(x) = \int_Ke^{(i\lambda-p)H(xk)}dk,\; \lambda
\in \mathfrak{a}^*_{\C},$$  $\rho =
\frac{1}{2}\sum_{\lambda\in\triangle^+} m_{\lambda}\cdot\lambda,$ where
$m_{\lambda}=dim (\mathfrak{g}_\lambda),$ and that $\varphi_{\lambda} =
\varphi_{\mu}$ iff $\lambda = s\mu$ for some $s \in \mathfrak{w}.$ Some of
the well-known properties of spherical functions are $\varphi_{-\lambda}(x^{-1}) =
\varphi_{\lambda}(x),$ $\varphi_{-\lambda}(x) =
\bar{\varphi}_{\bar{\lambda}}(x),$ $\mid \varphi_{\lambda}(x) \mid\leq \varphi_{\Re\lambda}(x),$ $\mid \varphi_{\lambda}(x)\mid\leq \varphi_{i\Im\lambda}(x),$ $\varphi_{-i\rho}(x)=1,$ $\lambda \in \mathfrak{a}^*_{\C},$ while $\mid \varphi_{\lambda}(x) \mid\leq \varphi_{0}(x),\;\lambda \in i\mathfrak{a}^{*},\;x \in G.$ Also if $\Omega$ is the \textit{Casimir operator} on $G$ then
$$\Omega\varphi_{\lambda} = -(\langle\lambda,\lambda\rangle +
\langle \rho, \rho\rangle)\varphi_{\lambda},$$ where $\lambda \in
\mathfrak{a}^*_{\C}$ and $\langle\lambda,\mu\rangle
:=tr(adH_{\lambda} \ adH_{\mu})$ for elements $H_{\lambda}$, $H_{\mu}
\in {\mathfrak{a}}.$ The elements $H_{\lambda}$, $H_{\mu}
\in {\mathfrak{a}}$  are uniquely defined by the requirement that $\lambda
(H)=tr(adH \ adH_{\lambda})$ and $\mu
(H)=tr(adH \ adH_{\mu})$ for every $H \in {\mathfrak{a}}$ ([$9.$],
Theorem $4.2$). Clearly $\Omega\varphi_0 = 0.$

Due to a hint dropped by Dixmier $[6.]$ $(cf.\;[12.])$ in his discussion of some functional calculus,
it is necessary to recall the notion of
a \textit{`positive-definite'} function and then discuss the situation for
positive-definite spherical functions.  We call a continuous function
$f: G \rightarrow \C$ (algebraically) positive-definite whenever, for all
$ x_1,\dots,x_m $ in $G$ and all $ \alpha_1,\dots,\alpha_m$ in $\C,$ we have $$\sum^m_{i,j=1}\alpha_i\bar{\alpha}_jf(x^{-1}_i x_j) \geq 0.$$  It
can be shown $(cf.\;[9.])$ that $f(e) \geq 0$ and $|f(x)| \leq f(e)$ for every
$x \in G$ implying that the space ${\wp}$ of all
positive-definite spherical functions on $G$ is a subset of the
space ${\mathfrak{F}}^{1}$ of all bounded spherical functions on $G.$

We know, by the Helgason-Johnson theorem ($[11.]$), that $${\mathfrak{F}}^{1}=
\mathfrak{a}^*+iC_{\rho}$$ where $C_{\rho}$ is the convex hull of $\{s\rho: s \in
\mathfrak{w}\}$ in $\mathfrak{a}^*.$ Defining the \textit{involution} $f^*$ of $f$ as $f^*(x) =
\overline{f(x^{-1})}$, it follows that $f = f^*$ for every $f \in
{\wp}$, and if $\varphi_{\lambda} \in {\wp}$, then $\lambda$
and $\bar{\lambda}$ are Weyl group conjugate, leading to a realization of $\wp$ as a subset of $\mathfrak{w} \setminus \mathfrak{a}^*_{\C}.$  ${\wp}$ becomes
a locally compact Hausdorff space when endowed with the \textit{weak $^{*}-$topology} as a subset of $L^{\infty}(G)$.

We denote the set of equivalence classes of the necessarily finite-dimensional irreducible representations of $K$ by $\mathcal{E}(K)$ whose character is $\chi_{\mathfrak{d}},$ for every $\mathfrak{d}\in\mathcal{E}(K).$ The class functions $\xi_{\mathfrak{d}}:K\rightarrow\C$ defined as $\xi_{\mathfrak{d}}(k):=dim(\mathfrak{d})\chi_{\mathfrak{d}}(k^{-1})$ are idempotents (i.e., $\xi_{\mathfrak{d}}\ast\xi_{\mathfrak{d}}=\xi_{\mathfrak{d}}$ with $\xi_{\mathfrak{d}_{1}}\ast\xi_{\mathfrak{d}_{2}}=0$ whenever $\mathfrak{d}_{1}\neq\mathfrak{d}_{2}$). Choosing $\pi$ to be any representation of $K$ (which may be the restriction to $K$ of a representation of $G$) in a complete Hausdorff locally convex space, $V,$ a continuous projection operator on $V$ may be given as the image of $\xi_{\mathfrak{d}}$ under $\pi.$ That is, $$E_{\pi,\mathfrak{d}}:=\pi(\xi_{\mathfrak{d}})=
\int_{K}\xi_{\mathfrak{d}}(k)\pi(k)dk=dim(\mathfrak{d})\int_{K}\chi_{\mathfrak{d}}(k^{-1})\pi(k)dk$$ (Here $\int_{K}dk=1$) Idempotency of $\xi_{\mathfrak{d}}$ assures that $E_{\pi,\mathfrak{d}}$ is indeed a projection on $V$ (since $E^{2}_{\pi,\mathfrak{d}}=E_{\pi,\mathfrak{d}}$ and $E_{\pi,\mathfrak{d}_{1}}E_{\pi,\mathfrak{d}_{2}}=0$ whenever $\mathfrak{d}_{1}\neq\mathfrak{d}_{2}$) and that its range, written as $V_{\mathfrak{d}}(=E_{\pi,\mathfrak{d}}(V))$ is a closed linear subspace of $V$ consisting mainly of members of $V$ which \textit{transform according to} $\mathfrak{d};\;[13.],\;p.\;109.$

The closed linear subspace $V_{\mathfrak{d}}$ above becomes familiar when $V=\mathcal{C}^{p}(G)$ under the usual regular representation. In this case the left and right regular representations are denoted as $l$ and $r$ given as $(l(x)f)(y)=f(x^{-1}y)$ and $(r(x)f)(y)=f(yx),$ respectively; $x,y\in G,\;f\in\mathcal{C}^{p}(G);$ and it may be computed that for any $\mathfrak{d}\in\mathcal{E}(K),$ $E_{l,\mathfrak{d}}=l(\xi_{\mathfrak{d}})$ and $E_{r,\mathfrak{d}}=r(\xi_{\mathfrak{d}})$ are the respective operators of left and right convolutions by the measure $\xi_{\mathfrak{d}}dk$ and $\overline{\xi}_{\mathfrak{d}}dk=\xi_{\overline{\mathfrak{d}}}dk,$ respectively. Here $\overline{\mathfrak{d}}$ is the class contragredient to $\mathfrak{d}.$ We therefore have a representation $l\times r$ of $G\times G$ on $\mathcal{C}(G)$ given as $((l\times r)(x,y)f)(z)=f(x^{-1}yz),$ $x,y,z\in G$ and the corresponding projection $E_{l\times r,(\mathfrak{d}_{1}\times\mathfrak{d}_{2})}=(l\times r)(\xi_{(\mathfrak{d}_{1}\times\mathfrak{d}_{2})}),$ which from the above remarks could be computed as $$E_{l\times r,(\mathfrak{d}_{1}\times\mathfrak{d}_{2})}f=\xi_{\mathfrak{d}_{1}}\ast f\ast\xi_{\mathfrak{d}_{2}},$$ $f\in\mathcal{C}^{p}(G).$

We now choose $\mathfrak{d}\in\mathcal{E}(K).$ The image of $\mathcal{C}^{p}(G)$ under $E_{l\times r,(\mathfrak{d}\times\overline{\mathfrak{d}})}$ is the closed subalgebra of $\mathcal{C}^{p}(G)$ denoted as $\mathcal{C}^{p}_{\mathfrak{d}}(G)$ and is exactly given as $$\mathcal{C}^{p}_{\mathfrak{d}}(G)=\xi_{\mathfrak{d}}\ast\mathcal{C}^{p}(G)\ast\xi_{\mathfrak{d}}=
\{\xi_{\mathfrak{d}}\ast f\ast\xi_{\mathfrak{d}}:\;f\in\mathcal{C}^{p}(G)\}.$$ Thus the members of $\mathcal{C}^{p}_{\mathfrak{d}}(G)$ are those of $\mathcal{C}^{p}(G)$ which may be written as $\xi_{\mathfrak{d}}\ast f\ast\xi_{\mathfrak{d}}$ for some $f\in\mathcal{C}^{p}(G).$ That is, every $g\in\mathcal{C}^{p}_{\mathfrak{d}}(G)$ is given as $g=\xi_{\mathfrak{d}}\ast f\ast\xi_{\mathfrak{d}},$  with $f\in\mathcal{C}^{p}(G).$ We shall henceforth write members $g$ of $\mathcal{C}^{p}_{\mathfrak{d}}(G)$ as $g_{\mathfrak{d},f}.$\\

{\bf $\bf{\S 3.}\;\;\;\;$ The Spherical Bochner Theorem and Extension}

Let $$\varphi_0(x):= \int_{K}\exp(-\rho(H(xk)))dk$$ be denoted
as $\Xi(x)$ and define $\sigma: G \rightarrow \C$ as
$$\sigma(x) = \|X\|$$ for every $x = k\exp X \in G,\;\; k \in K,\; X
\in \mathfrak{a},$ where $\|\cdot\|$ is a norm on the finite-dimensional
space $\mathfrak{a}.$ These two functions are spherical functions on
$G$ and there exist numbers $c,d$ such that $$1 \leq \Xi(a)
e^{\rho(\log a)} \leq c(1+\sigma(a))^d.$$ Also there exists $r_0
> 0$ such that $c_0 =: \int_G\Xi(x)^2(1+\sigma(x))^{r_0}dx
< \infty$ ($[16.],$ p. $231$).  For each
$0 \leq p \leq 2$ define ${\cal C}^p(G)$ to be the set consisting of
functions $f$ in $C^{\infty}(G)$ for which $$\|f\|_{g_1,
g_2;m} :=\sup_G|f(g_1; x ; g_2)|\Xi (x)^{-2/p}(1+\sigma(x))^m <
\infty$$ where $g_1,g_2 \in \mathfrak{U}(\mathfrak{g}_{\C}),$ the \textit{universal
enveloping algebra} of $\mathfrak{g}_{\C},$ $m \in \Z^+, x \in G,$
$f(x;g_2) := \left.\frac{d}{dt}\right|_{t=0}f(x\cdot(\exp tg_2))$
and $f(g_1;x) :=\left.\frac{d}{dt}\right|_{t=0}f((\exp
tg_1)\cdot x).$ We call ${\cal C}^p(G)$ the Schwartz space on $G$
for each $0 < p \leq 2$ and note that ${\cal C}^2(G)$ is the
well-known (see $[1.]\;\mbox{and}\;[2.]$) Harish-Chandra space of rapidly decreasing functions on
$G.$ The inclusions $$C^{\infty}_{c}(G) \subset {\cal C}^p(G)
\subset L^p(G)$$ hold and with dense images. It also follows that
${\cal C}^p(G) \subseteq {\cal C}^q(G)$ whenever $0 \leq p \leq q
\leq 2.$ Each ${\cal C}^p(G)$ is closed under \textit{involution} and the
\textit{convolution}, $*.$ Indeed ${\cal C}^p(G)$ is a Fr$\acute{e}$chet algebra ($[15.],$ p. $69$). We endow ${\cal C}^p(G//K)$
with the relative topology as a subset of ${\cal C}^p(G)$.

For any measurable function $f$ on $G$ we define the \textit{Harish-Chandra Fourier
transform} $\widehat{f}$ as $$\widehat{f}(\lambda) = \int_G f(x)
\varphi_{-\lambda}(x)dx,$$ $\lambda \in \mathfrak{a}^*_{\C}.$ We shall call it \textit{spherical} whenever $f$ is $K-$biinvariant. It is known (see $[3.]$) that for $f,g \in L^1(G)$ we have:\\
\begin{enumerate}
\item [$(i.)$] $(f*g)^{\wedge} = \widehat{f}\cdot\widehat{g}$ on $ {\mathfrak{F}}^{1}$
whenever $f$ (or $g$) is right - (or left-) $K$-invariant; \item
[$(ii.)$] $(f^*)^{\wedge}(\varphi) =
\overline{\widehat{f}(\varphi^*)}, \varphi \in {\mathfrak{F}}^{1}$; hence
$(f^*)^{\wedge} = \overline{\widehat{f}}$ on ${\wp}:$ and, if we
define $$f^{\#}(g) := \int_{K\times K}f(k_1xk_2)dk_1dk_2,  x\in
G,$$ then \item [$(iii.)$] $(f^{\#})^{\wedge}=\widehat{f}$ on ${\mathfrak{F}}^{1}.$
\end{enumerate}

In order to know the image of the spherical Fourier transform when
restricted to ${\cal C}^p(G//K)$ we need the following spaces that are central to the statement
of the well-known result of Trombi and Varadarajan [$14.$] (Theorem $3.1$ below).

Let $C_\rho$ be the closed convex hull of the (finite) set $\{s\rho :
s\in \mathfrak{w}\}$ in $\mathfrak{a}^*$, i.e., $$C_\rho =
\left\{\sum^n_{i=1}\lambda_i(s_i\rho) : \lambda_i \geq 0,\;\;\sum^n_{i=1}\lambda_i = 1,\;\;s_i \in \mathfrak{w}\right\}$$ where we recall that, for every
$H \in \mathfrak{a},$ $$(s\rho)(H) = \frac{1}{2} \sum_{\lambda\in\triangle^+}
 m_{\lambda}\cdot\lambda (s^{-1}H).$$  Now for each
$\epsilon > 0$ set ${\mathfrak{F}}^{\epsilon} = \mathfrak{a}^*+i\epsilon
C_\rho.$ Each ${\mathfrak{F}}^{\epsilon}$ is convex in $\mathfrak{a}^*_{\C}$ and
$$int({\mathfrak{F}}^{\epsilon}) =
\bigcup_{0<\epsilon'<\epsilon}{\mathfrak{F}}^{\epsilon^{'}}$$
([$14.$], Lemma $(3.2.2)$). Let us define $\mathcal{Z}({\mathfrak{F}}^{0}) = \mathcal{S}
(\mathfrak{a}^*)$ and, for each $\epsilon>0,$ let
$\mathcal{Z}({\mathfrak{F}}^{\epsilon})$ be the space of all $\C$-valued
functions $\Phi$ such that  $(i.)$ $\Phi$ is defined and holomorphic
on $int({\mathfrak{F}}^{\epsilon}),$ and $(ii.)$ for each holomorphic
differential operator $D$ with polynomial coefficients we have $\sup_{int({\mathfrak{F}}^{\epsilon})}|D\Phi| < \infty.$ The space
$\mathcal{Z}({\mathfrak{F}}^{\epsilon})$ is converted to a Fr$\acute{e}$chet algebra by equipping it with the
topology generated by the collection, $\| \cdot \|_{\mathcal{Z}({\mathfrak{F}}^{\epsilon})},$ of seminorms given by $\|\Phi\|_{\mathcal{Z}({\mathfrak{F}}^{\epsilon})} := \sup_{int({\mathfrak{F}}^{\epsilon})}|D\Phi|.$  It is known that $D\Phi$ above extends to
a continuous function on all of ${\mathfrak{F}}^{\epsilon}$
([$14.$], pp. $278-279$).  An appropriate subalgebra of
$\mathcal{Z}({\mathfrak{F}}^{\epsilon})$ for our purpose is the closed
subalgebra $\bar{\mathcal{Z}}({\mathfrak{F}}^{\epsilon})$ consisting of
$\mathfrak{w}$-invariant elements of $\mathcal{Z}({\mathfrak{F}}^{\epsilon})$,
$\epsilon \geq 0.$ The following well-known result affords us the
opportunity of defining a distribution on ${\cal C}^p(G//K).$

{\bf 3.1 Theorem (Trombi-Varadarajan $[14.]$).}  \textit{Let $0 < p \leq 2$ and
set $\epsilon = (\frac{2}{p})-1$.  Then the spherical Harish-Chandra
Fourier transform $f \mapsto \widehat{f}$ is a linear
topological algebra isomorphism of ${\cal C}^p(G//K)$ onto $\bar{\mathcal{Z}}
({\mathfrak{F}}^{\epsilon}).\;\;\Box$}

It has been shown in $[13.]$ that the function $(\widehat{\xi_{\mathfrak{d}}})^{-1}$ exists and lives on $K/M\times int(\mathfrak{F}^{\epsilon})$ and if we define $\mathcal{C}^{p}(\widehat{G})$ as $$\mathcal{C}^{p}(\widehat{G})=
\{(\widehat{\xi_{1}})^{-1}\cdot h\cdot(\widehat{\xi_{1}})^{-1}:\;h\in\bar{\mathcal{Z}}({\mathfrak{F}}^{\epsilon})\}$$ for $0<p\leq 2,$ then a complete generalization of Theorem $3.1$ to all of ${\cal C}^p(G)$ is possible, allowing us define a distribution on ${\cal C}^p(G),$ and is given as follows.

\textbf{3.2. Theorem (The Fundamental Theorem of Harmonic Analysis on $G$ [13.]).} \textit{Let $0<p\leq 2,$ then the Harish-Chandra Fourier transform sets up a linear topological algebra isomorphism $\mathcal{C}^{p}(G)\rightarrow\mathcal{C}^{p}(\widehat{G}).\;\Box$}

In order to use the above theorems to re-state and generalize the results of Barker $[3.],$ we require the following notions.

{\bf 3.3 Definitions.}
\begin{enumerate}
\item [$(i.)$] \textit{A distribution  $S$ on $G$ (i.e., $S \in C^{\infty}_c(G)'$) is said to be (integrally) \textit{positive-definite} (written as $S\gg 0$) whenever $$S [f * f^*] \geq 0,$$ for $f \in C^{\infty}_c(G).$}
    \item [$(ii.)$] \textit{A distribution $S$ on $G$ is called $K$-\textit{bi-invariant} whenever $S^{\#}=S$ where
$$S^{\#}[f] := S[f^{L(k_1)R(k_2)}],$$ for $f \in C^{\infty}_c(G).$}
\item [$(iii.)$] \textit{A measure $\mu$ defined on ${\wp}$ is said to be of \textit{polynomial growth} if there exists a holomorphic polynomial $Q$ on $\mathfrak{a}^*_{\C}$ such that $\int_{\wp}(d\mu/|Q|)<\infty.$}
    \item [$(iv.)$] \textit{The \textit{support}, $supp (\mu),$ of a regular Borel measure $\mu$ is the smallest closed set $A$ such that $\mu (B)= 0$ for all Borel sets $B$ disjoint from $A.$}
    \end{enumerate}

One of the two well-known most general results on positive-definite distributions on $G$ (indeed, on any unimodular Lie group) and which started off investigations leading to Bochner theorem (on $\R^{n}$) is the following generalization by Barker $[4.]$ used in the proof of his Bochner theorems.

\textbf{3.4 Theorem $[4.].$} Every $S\in (C^{\infty}_c(G))'$ in which $S \gg 0$ can be expressed as the finite sum of $$S=\sum_{j}D^{j}E^{j}f_{j}$$ where each $f_{j}$ is a bounded function on $G$ and $D^{j}$ (respectively, $E^{j}$) is a left (respectively, right) invariant differential operators$.\;\Box$

The following is the first of the main results of $[3.]$ giving the spherical Bochner theorem ($cf.\;[8.]$).

{\bf 3.5 Theorem (The spherical Bochner theorem) $[3.].$}  \textit{Let $S \in (C^{\infty}_c(G))'$ in which $S \gg 0$.
Then $S$ extends uniquely to an element in $(\mathcal{C}^1(G))'$ and there exists a unique positive regular Borel measure $\mu$ of polynomial growth on ${\wp}$ such that $$S[f] = \int_{\wp}\widehat{f}d\mu,\;\; f \in {\cal C}^1(G//K).$$ The correspondence between $S$ and $\mu$ is bijective when restricted to $K-$bi-invariant distributions, in which case the formula holds for all $f \in {\cal C}^1(G).\;\;\Box$}

The unique Borel measure $\mu$ in the last theorem shall be called the (spherical) Bochner measure of $T.$ The second of the main results of $[3.]$ (which relates the largest $\mathcal{C}^{p}(G//K)-$space to which the positive-definite distribution $S$ of Theorem $3.4$ can be extended with the support of its (spherical) Bochner measure) is also a consequence of the Trombi-Varadarajan theorem (Theorem $3.1$ above) and is stated as follows.

{\bf 3.6 Theorem (The extension theorem) $[3.].$}  \textit{Suppose $S$ is a positive-definite distribution with spherical Bochner measure $\mu$.  Then $S \in ({\cal C}^p(G//K))'$ if and only if supp $(\mu) \subset {\mathfrak{F}}^{\epsilon}$ where $1 \leq p \leq 2$ and $\epsilon =(\frac{2}{p})-1$.  In such a case $$S[f] = \int_{\wp}\widehat{f}d\mu,\;\;f \in {\cal C}^p(G//K).\;\;\Box$$}

It is our modest aim in the next section to show that Barker $[3.]$ had the needed framework to arrive at the full Bochner theorem on $G,$ but was hampered by the non-availability of the full theory for the Harish-Chandra Fourier transform on $\mathcal{C}^{p}(G)$ (given here as Theorem $3.2$) and was therefore restricted to the spherical Bochner theorem derived through Theorem $3.1.$ It is now possible, in the light of Theorem $3.2$ to investigate the full Bochner theorem as done in the next section.\\

\textbf{\S 4. The full Bochner theorem.}

We shall in the present section, establish the full Bochner theorem effortlessly from the results of \S 3. It will be clear from our results that the spherical Bochner measure corresponding to a positive-definite distribution on $G$ (of Theorems $3.5$ and $3.6$ above) is still the positive regular Borel measure for each $T\in(\mathcal{C}^p(G))'.$ An indication of this fact is however already contained in the first part of Theorem $3.5.$ Even though the (so-called \textit{spherical}) Bochner measure is still defined in this section on $\wp$ (indeed, supp $(\mu)=\wp$), we shall henceforth refer to it simply as \textit{the Bochner measure}. This is partly due to the main result of $[13]$ (which is Theorem $3.2$ above) which shows that the Trombi-Varadarajan spherical image $\bar{\mathcal{Z}}
({\mathfrak{F}}^{\epsilon})$ is indeed the \textit{nucleus} of the image of the full Harish-Chandra Fourier transform on all of $G$ and partly due to the next theorem.

We start by recalling a continuous surjection map of $\mathcal{C}^p(G)$ onto $\mathcal{C}^p(G//K).$

\textbf{4.1 Lemma.} \textit{Let $\varphi\in C^{\infty}_{c}(G)$ and consider $\varphi^{\#}\in C^{\infty}_{c}(G//K).$ Then $\varphi\mapsto\varphi^{\#}$ is a continuous linear map of $\mathcal{C}^p(G)$ onto $\mathcal{C}^p(G//K).$}

\textbf{Proof.} $\varphi\mapsto\varphi^{\#}$ is a continuous linear map of $C^{\infty}_{c}(G)$ onto $C^{\infty}_{c}(G//K),$ which extends to $\mathcal{C}^p(G)$ due to the denseness of $C^{\infty}_{c}(G)$ in $\mathcal{C}^p(G).\;\Box$

Indeed, $\varphi\mapsto\varphi^{\#}$ is a continuous endomorphism of $\mathcal{C}^p(G)$ onto $\mathcal{C}^p(G//K),$ as proved in Proposition $3.2$ of $[3.].$ The following is our first main result.

\textbf{4.2 Theorem (The full Bochner theorem.)} \textit{Let $T \in (C^{\infty}_c(G))'$ in which $T \gg 0$.
Then $T$ extends uniquely to an element in $(\mathcal{C}^1(G))'$ and there exists a unique positive regular Borel measure $\mu$ of polynomial growth on ${\wp}$ such that $$T[f] = \int_{\wp}\widehat{f}d\mu,\;\; f \in {\cal C}^1(G).$$ The correspondence between $T$ and $\mu$ is bijective when restricted to $K-$bi-invariant distributions.}

\textbf{Proof.} Assertion $(i.)$ follows from the first part of Theorem $3.5.$ For $(ii.)$ let $f\in{\cal C}^1(G).$ Then by Lemma $4.1,$ $f\mapsto f^{\#}$ is a continuous linear map of $\mathcal{C}^1(G)$ onto $\mathcal{C}^1(G//K)$ and $\widehat{(f^{\#})}=\widehat{f}$ on $\mathfrak{F}^{1},$ hence on $\wp$ (by Helgason-Johnson theorem). This means that the spherical Harish-Chandra Fourier transform of the spherical function $f^{\#}$ is exactly the full Harish-Chandra Fourier transform of the function $f$ on $G.$

Now, since $T\in(\mathcal{C}^1(G))',$ it follows from Proposition $3.2$ of $[3.]$ that $T_{|_{\mathcal{C}^1(G//K)}}\in(\mathcal{C}^1(G//K))'.$ We shall denote $T_{|_{\mathcal{C}^1(G//K)}}$ by $S_{T}$ and note that it is such that $T[f]=S_{T}[f^{\#}]$ for all $f\in{\cal C}^1(G).$ It follows from Theorem $3.5$ therefore that $$T[f]=S_{T}[f^{\#}]=\int_{\wp}\widehat{(f^{\#})}d\mu=\int_{\wp}\widehat{f}d\mu,$$ as required for all $f\in{\cal C}^1(G).$ The claim in $(iii.)$ holds by viewing the bijective correspondence of Theorem $3.5$ in the light of $(ii.).\;\Box$

The unique positive regular Borel measure $\mu$ in Theorem $4.2$ shall simply be called the Bochner measure of $T.$ In the light of Theorem $4.2$ we shall now prove the full version of Barker's Extension (Theorem $3.6$) of the spherical Bochner theorem as follows. This will be our first extension of Theorem $4.2.$

\textbf{4.3 Theorem (First extension theorem.)} \textit{Suppose $T$ is a positive-definite distribution with Bochner measure $\mu$.  Then $T \in ({\cal C}^p(G))'$ if and only if supp $(\mu) \subset {\mathfrak{F}}^{\epsilon}$ where $1 \leq p \leq 2$ and $\epsilon = (\frac{2}{p})-1$.  In this case $$T[f]=\int_{\wp}\widehat{f}d\mu,\;\;f \in {\cal C}^p(G).$$}

\textbf{Proof.} Given that $T \in ({\cal C}^p(G))'$ satisfies $T\gg0,$ then $S_{T}$ defined as $S_{T}:=T_{|_{\mathcal{C}^p(G//K)}}$ satisfies $S_{T}\in(\mathcal{C}^p(G//K))'$ with $S\gg0,$ corresponding to the (spherical) Bochner measure $\mu.$ It follows therefore from Theorem $3.6$ that supp $(\mu) \subset {\mathfrak{F}}^{\epsilon}.$

Conversely, let $\mu$ denote the Bochner measure of a positive-definite distribution $S$ on $G$ with supp $(\mu) \subset {\mathfrak{F}}^{\epsilon}.$ The by Theorem $3.6,$ $S\in(\mathcal{C}^p(G//K))'$ and $S[f] = \int_{\wp}\widehat{f}d\mu,$ for all $f \in {\cal C}^p(G//K).$ If we now define $T$ as $T[f]:=S[f^{\#}]$ for all $f \in {\cal C}^p(G),$ then; $(i.)$ $T$ is a distribution on $G;$ $(ii.)$ $T$ is positive-definite with Bochner measure $\mu,$ since $T[f\ast f^{\ast}]=S[(f\ast f^{\ast})^{\#}]=\int_{\wp}\widehat{(f\ast f^{\ast})^{\#}}d\mu=\int_{\wp}\widehat{(f\ast f^{\ast})}d\mu=\int_{\wp}\mid\widehat{f}\mid^{2}d\mu\geq0;$ and $(iii.)$ for all $f \in {\cal C}^p(G),$ $T[f]:=S[f^{\#}]=\int_{\wp}\widehat{(f^{\#})}d\mu=\int_{\wp}\widehat{f}d\mu.\;\Box$

\textbf{4.4 Remarks.} $(i.)$ With $p=2$ Theorem $4.3$ establishes Theorem $9.3$ of $[5.]$ for $G=SU(1,1)/\{\pm1\}.$

$(ii.)$ The above theorem gives a characterization of any given $p-$tempered positive-definite distribution in terms of which of the tubular domains, $\mathfrak{F}^{\epsilon},$ contains the support of its corresponding Bochner measure, $\mu.$ The existence and explicit realization for $p=1$ is given in Theorem $4.2,$ where supp $(\mu)=\wp\subset\mathfrak{F}^{1}.$ Here and by the Helgason-Johnson theorem, $\mathfrak{F}^{1}$ is the exclusive indexing set for all \textit{bounded} spherical functions, in terms of which every positive-definite distribution on $G$ is written, according to Theorem $3.4.$ It is therefore safe to say that, unless Theorem $3.4$ is generalized, the first extension theorem (Theorem $4.3$) above may not be fully realized, as (in any analysis based on Theorem $3.4$) we shall always be taken back to the nature of the functions $f_{j}$ as being bounded functions; thus $\mathfrak{F}^{1}$ may equally surface again and leading eventually to $\mathcal{C}^1(G).$

Indeed any change in the nature of $f_{j}$ (in Theorem $3.4$) is bound to reflect in the type of Bochner theorem (and extension) that would be available for consideration. A general look at the restriction introduced through Theorem $3.4$ (though a very general result in its own right) is expected to follow the following line of thought.

\textbf{4.5 Conjecture (Abstract extension theorem).} \textit{Let $T$ denote a positive-definite distribution on $G$ with Bochner measure $\mu$ and set $\epsilon=(\frac{2}{p})-1$ for all $1\leq p\leq2.$ Then $T \in ({\cal C}^p(G))'$ if and only if $T$ can be expressed as a finite sum of derivatives of functions that are exclusively indexed by $\mathfrak{F}^{\epsilon}$ with supp $(\mu) \subset {\mathfrak{F}}^{\epsilon}.\;\Box$}

The proof and application of this conjecture may be useful in rectifying the anomaly of the first parts of Theorems $3.5$ and $4.2,$ by showing that every positive-definite distribution is $2-$tempered (as known for $\R^{n}$). That is, that the support of the Bochner measure is always in ${\mathfrak{F}}^{0}.$

The following is a consequence of Theorem $4.3$ and the fundamental theorem (Theorem $3.2$).

\textbf{4.6 Corollary.} \textit{Let $T\in(\mathcal{C}^{p}(G))',\;T\gg0,$ with Bochner measure $\mu$ and set $\epsilon=(\frac{2}{p})-1$ for all $1\leq p\leq2$ such that supp $(\mu) \subset {\mathfrak{F}}^{\epsilon}.$ Then, for all $f\in\mathcal{C}^{p}(G),$ $$T[f]=\int_{\wp}(\widehat{\xi_{1}})^{-1}\cdot\widehat{\varphi}\cdot(\widehat{\xi_{1}})^{-1}d\mu,$$ where $\varphi=f^{\#}.\;\Box.$}

We may now embark on an application of Theorem $4.3$ to $p-$tempered invariant positive-definite distributions in the fashion of $[5.].$ This is possible now that, according to Theorem $4.3,$ every $p-$tempered positive-definite distribution on $G$ is the (full) Harish-Chandra Fourier transform of some unique tempered measure $\mu$ (of course, with supp $(\mu) \subset {\mathfrak{F}}^{\epsilon}$). With the (full) Harish-Chandra Fourier transform, it will be possible to then bring the full Representation Theory of $G$ to bear on more explicit expression for and decomposition of $T[f],$ for all $f\in\mathcal{C}^{p}(G),$ as is known in $[5.].$\\

{\bf \S 5. References.}
\begin{description}
\item [{[1.]}] Arthur, J. G., \textit{Harmonic analysis of tempered distributions on semisimple Lie groups of real rank one,} Ph.D. Dissertation, Yale University, $1970.$
    \item [{[2.]}] Arthur, J. G., \textit{Harmonic analysis of the Schwartz space of a reductive Lie group,} I. II. (preprint, $1973$).
        \item [{[3.]}] Barker, W. H.,  The spherical Bochner theorem on semisimple Lie groups, \textit{J. Funct. Anal.,} vol. \textbf{20}  ($1975$), pp. $179-207.$
            \item [{[4.]}] Barker, W. H.,  Positive definite distributions on unimodular Lie groups, \textit{Duke Mathematical Journal,} vol. \textbf{43} No. \textbf{1},, ($1976$), pp. $71-79.$
            \item [{[5.]}] Barker, W. H.,  Tempered invariant, positive-definite distributions on\\ $SU(1,1)/\{\pm1\},$
            \textit{Illinois J. Maths}, vol. \textbf{28}, no. $1,$ ($1984$), pp. $83-102.$
                \item [{[6.]}] Dixmier, J., Op$\acute{\mbox{e}}$rateurs de rang fini dans les repr$\acute{\mbox{e}}$sentations unitaires,\textit{ Publ. math. de l' Inst. Hautes $\acute{\mbox{E}}$tudes Scient.,} tome $\textbf{6}$ ($1960$), pp. $13-25.$
                    \item [{[7.]}] Godement, R. A.,  A theory of spherical functions, $I.$ \textit{Trans. Amer. Math. Soc.,} vol. \textbf{73} ($1952$), pp. $496-556.$
                        \item [{[8.]}] Godement, R. A., Introduction aux travaux de A. Selberg, \textit{S$\acute{e}$minaire Bourbaki} No. \textbf{144} Paris ($1957$).
                    \item [{[9.]}] Helgason, S., \textit{``Differential Geometry, Lie Groups and Symmetric Spaces,''} Academic Press, New York, $1978.$
                \item [{[10.]}] Helgason, S., \textit{``Groups and Geometric Analysis; Integral Geometry, Invariant Differential Operators, and Spherical Functions,''} Academic Press, New York, $1984.$
                \item [{[11.]}] Helgason, S. and Johnson, K.,  The bounded spherical functions on symmetric spaces,  \textit{Advances in Math,} \textbf{3} ($1969$), pp. 586-593.
                        \item [{[12.]}] Kahane J.-P.,  Sur un th$\acute{e}$or$\grave{e}$me de Wiener-L$\acute{e}$vy, \textit{C.R. Acad. Sc., Paris,} t. $\textbf{246},$ ($1958$), pp. $1949-1951.$
                            \item [{[13.]}] Oyadare, O. O.,  Non-spherical Harish-Chandra Fourier transforms on real reductive groups, arXiv:$1907.00717$v3 [mathFA] $16$ Jul 2019. Submitted for review.
                    \item [{[14.]}] Trombi, P. C. and Varadarajan, V. S.,  Spherical transforms on semisimple Lie groups, \textit{Ann. of Math.,} \textbf{94} ($1971$), pp. $246-303.$
                \item [{[15.]}] Varadarajan, V. S.,  The theory of characters and the discrete series for semisimple Lie groups, in \textit{Harmonic Analysis on Homogeneous Spaces,} (C.C.  Moore (ed.)) \textit{Proc. of Symposia in Pure Maths.,} vol. \textbf{26} ($1973$), pp. $45-99.$
                    \item [{[16.]}] Varadarajan, V. S., \textit{``An introduction to harmonic analysis on semisimple Lie groups,''}  Cambridge University Press, Cambridge, $1989.$
                        \end{description}

\end{document}